\documentstyle{amsart}

\title[Anisotropic discrete problems] {On sequences of solutions for discrete\\ anisotropic equations
}

\author{Giovanni Molica Bisci}
\address[G. Molica Bisci]{Dipartimento P.A.U., Universit\`a  degli
Studi Mediterranea di Reggio Calabria, Salita Melissari - Feo di
Vito, 89100 Reggio Calabria, Italy} \email{gmolica@@unirc.it}
\author{Du\v{s}an Repov\v{s}}
\address[D. Repov\v{s}]{Faculty of of Education, and Faculty of Mathematics and Physics, University of Ljubljana, Kardeljeva pl. 16, Ljubljana, Slovenia 1000
}
\email{dusan.repovs@@guest.arnes.si}

\thanks{{\em 2010 Mathematics Subject Classification:} 47A75, 35B38, 35P30, 34L05, 34L30.}

\keywords{Discrete nonlinear boundary value problems; anisotropic difference equations;
 infinitely
many solutions; variational methods.}

{}
\newtheorem{theorem}{Theorem}[section]
\newtheorem{lemma}{Lemma}[section]
{}

\newtheorem{corollary}{Corollary}[section]
\newtheorem{remark}{Remark}[section]
\newtheorem{example}{Example}[section]

\newcommand{\RR}{\mbox{\normalshape I\!R}}
\newcommand{\erre}{\mbox{\normalshape I\!R}}
\newcommand{\enne}{\mbox{\normalshape I\!N}}
\def\ZZ{{\mathbb Z} }

\begin{document}

\begin{abstract}
Taking advantage of a recent critical point theorem, the existence of infinitely many solutions for an anisotropic problem with a parameter is established.
More precisely, a concrete interval of positive parameters, for which the treated problem admits infinitely many solutions, is determined without symmetry assumptions on the nonlinear data. Our goal was achieved by requiring an appropriate behavior of the nonlinear terms at zero, without any additional conditions.
\end{abstract}
\maketitle
\section{Introduction}
For every $a,b\in \ZZ$, such that $a<b$, set
$\ZZ[a,b]:=\{a,a+1,...,b\}$.
This work is concerned with the study of existence of solutions for the following anisotropic difference equation
\begin{equation}\tag{$A_{\lambda}^{f}$} \label{ani}
\left\{\begin{array}{lll}
-\Delta(|\Delta u(k-1)|^{p(k-1)-2}\Delta u(k-1))=\lambda f_k(u(k)),\,\,
k\in\ZZ[1,T]\\
u(0)=u(T+1)=0,
\end{array}\right.
\end{equation}
where $\lambda$ is a positive parameter, $f_k:\erre\rightarrow \erre$ is a continuous function for every $k\in\ZZ[1,T]$ (with $T\geq 2$), and
$\Delta u(k-1):=u(k)-u(k-1)$ is the
forward difference operator. Further, we will assume that the map
$p:\ZZ[0,T+1]\rightarrow \erre$ satisfies $p^-:=\displaystyle\min_{\ZZ[0,T+1]}p(k)>1$ as well as $p^+:=\displaystyle\max_{\ZZ[0,T+1]}p(k)>1$.\par
Discrete boundary value problems have been intensively studied in
the last decade. The modeling of certain nonlinear problems from
biological neural networks, economics, optimal control and other
areas of study have led to the rapid development
 of the theory of difference equations; see the monograph of Agarwal \cite{A}.
\par
Our idea here is to transfer
the problem of existence of solutions for problem \eqref{ani} into
the problem of existence of critical points for a suitable associated
energy functional, namely $J_\lambda$.
More precisely, the main purpose of this paper is to investigate the existence of infinitely many solutions to problem \eqref{ani} by using a critical point theorem obtained in \cite{Ricceri}; see Theorem \ref{two}.\par
Continuous analogues of problems like \eqref{ani} are known to be mathematical models of various phenomena arising in the study
of elastic mechanics, electrorheological fluids or image restoration; see, for instance, Zhikov \cite{Z}, R\.{u}\v{z}i\v{c}ka \cite{R} and Chen, Levine and Rao \cite{CLR}.\par Variational continuous
anisotropic problems have been studied by Fan and Zhang in \cite{FZ} and later considered by many methods and authors, see \cite{HH}
for an extensive survey of such boundary value problems.\par
 Research concerning the discrete anisotropic problems of
type \eqref{ani} was initiated by Kone and Ouaro in \cite{KO} and by Mih\u{a}ilescu, R\u{a}dulescu and Tersian in \cite{MRT}. In these papers known tools from the critical point theory are applied in order to get the
existence of solutions.\par

We also recall that very recently, in \cite{GG}, Galewski and G{\l}\c{a}b considered the discrete anisotropic boundary value problem \eqref{ani} using critical
point theory. First,
 they applied the direct method of the calculus of variations and the
mountain pass technique in order to reach the existence of at least one nontrivial solution.
Moreover, in the same paper, a discrete three critical point theorem was exploited
in order to get the existence of at least two nontrivial solutions.\par

\indent We note that most existence results for discrete problems assume that the nonlinearities data are odd functions. Only a few papers deal with nonlinearities for which this property does not hold; see, for instance, the interesting paper of Krist\'{a}ly, Mih\u{a}ilescu and R\u{a}dulescu \cite{KMR}; see also \cite{KMRT}.\par
 In analogy with the cited contributions, in our approach we do not require any symmetry hypothesis. A special case of our contributions reads as follows.

\begin{theorem}\label{intro}Let $g:\erre\rightarrow \erre$ be a nonnegative and continuous function. Assume that
\[
\displaystyle{\liminf_{t\rightarrow 0^+}}\frac{\displaystyle\int_{0}^{t}g(s)ds}{t^{p^+}}=0\,\,\,\,\,\,{\rm and}\,\,\,\,\,\limsup_{t\rightarrow 0^+}\frac{\displaystyle\int_{0}^{t}g(s)ds}{t^{p^-}}=+\infty.
\]
\noindent Then, for each $\lambda>0$, the problem
\begin{equation}\tag{$A_{\lambda}^{g}$} \label{ani2}
\left\{\begin{array}{lll}
-\Delta(|\Delta u(k-1)|^{p(k-1)-2}\Delta u(k-1))=\lambda g(u(k)),\,\,
k\in\ZZ[1,T]\\
u(0)=u(T+1)=0,
\end{array}\right.
\end{equation}
admits a sequence of nonzero solutions which converges to zero.
\end{theorem}
The plan of the paper is the following: Section 2 is devoted to our abstract framework, while Section 3 is dedicated to main results. A concrete example of application of the attained abstract results is presented; see Example \ref{esempio}.\par
 Finally, we cite the recent monograph by Krist\'aly, R\u adulescu and Varga \cite{KRV} as a
general reference on variational methods adopted here.

\section{Abstract framework}

Let $(X,\|\cdot\|)$ be a finite dimensional Banach space and let
$J_\lambda:X\rightarrow \erre$ be a function satisfying the
following structure hypothesis:
\begin{itemize}
\item[$(\Lambda)$]\textit{$J_\lambda(u):=\Phi(u)-\lambda\Psi(u)$
for all $u\in X$, where $\Phi, \Psi:X\rightarrow \erre$ are
functions of class $C^1$ on $X$ with $\Phi$ coercive, i.e.
$\lim_{\|u\|\rightarrow \infty}\Phi(u)=+\infty$, and $\lambda$ is
a real positive parameter.}
\end{itemize}
\noindent Moreover, provided that $r>\inf_{X}\Phi$, put
\[
\varphi(r):=\inf_{u \in \Phi^{-1}(\left]-\infty, r\right[)}
\frac{\left(\displaystyle\sup_{v \in {\Phi^{-1}\left(\left]-\infty,
r \right[ \right)}}\Psi(v)\right)-\Psi(u)}{r - \Phi(u)},
\]
\noindent and
$$
\delta:=\liminf_{r\rightarrow (\inf_X\Phi)^{+}}\varphi(r).
$$
\indent Clearly, one can observe that $\delta \ge 0$. Further, when $\delta = 0$, in the sequel, we agree that
$1/\delta$ is $+\infty$.\par
\begin{theorem} \label{two}
Assume that the condition $(\Lambda)$ is verified. If $\delta<+\infty$ then for each
$\lambda\in \left]0,{1}/{\delta}\right[$, one of the following
holds$:$\\
either
\begin{itemize}
 \item[$(\textrm{b}_1)$] there is a global minimum of $\Phi$ which is a local minimum of $J_\lambda$,
\end{itemize}
or
\begin{itemize}
\item[$(\textrm{b}_2)$]there is a sequence $\{u_m\}$ of pairwise distinct critical points
  $($local minima$)$ of $J_\lambda$, with $\lim_{m \to \infty}\Phi(u_m)=\inf_X\Phi$, which
  converges to a global minimum of $\Phi$.
\end{itemize}
\end{theorem}
\begin{remark}\rm{Theorem \ref{two}
is a special form of the quoted variational principle of Ricceri contained in \cite{Ricceri}.}
\end{remark}

By a \textit{strong solution} (briefly called a ``solution") to \eqref{ani} we mean such a function $u:\ZZ[0,T+1]\rightarrow \erre$
which satisfies the given equation and the associated
boundary conditions.
{Solutions will be investigated in the space} $$H=\{u:\ZZ[0,T+1]\rightarrow\RR;\
u(0)=u(T+1)=0\} .$$

\indent Clearly, $H$ is a $T$-dimensional Hilbert
space with the inner product
$$(u,v):=\sum_{k=1}^{T+1}\Delta u(k-1)\Delta
v(k-1),\;\;\;\forall\;u,v\in H.$$ The associated norm is defined
by $$\|u\|=\left(\sum_{k=1}^{T+1}|\Delta
u(k-1)|^2\right)^{1/2}\,,$$
see the work \cite{APR} for details.\par
\indent From now on, for every $u\in H$, set
\begin{equation*}
\Phi(u):=\sum_{k=1}^{T+1}\frac{1}{p(k-1)}|\Delta
u(k-1)|^{p(k-1)},\,\,\,\,\,\mbox{and}\,\,\,\,\Psi(u):=\sum_{k=1}^{T}F_k(u(k)),
\end{equation*}
where
$F_k(t):=\displaystyle\int_{0}^{t}f_k(s)ds$, for every $t\in \erre$ and $k\in \ZZ[1,T]$.
Further, let us denote
 $$
J_\lambda(u):=\Phi(u)-\lambda\Psi(u),
$$
for every $u\in H$.\par
 \indent Standard arguments assure that $J_\lambda\in C^1(H;\RR)$ and
$$\langle J'_\lambda(u),v\rangle=\sum_{k=1}^{T+1}|\Delta u(k-1)|^{p(k-1)-2}\Delta u(k-1)\Delta
v(k-1)-\lambda\sum_{k=1}^{T}
f_k(u(k))v(k),$$
for all $u,v\in H.$\par
Next result will be crucial in the sequel.

\begin{lemma}\label{LL} The functional $\Phi:H\rightarrow \erre$ is coercive, i.e.
$$
\lim_{\|u\|\rightarrow \infty}\sum_{k=1}^{T+1}\frac{1}{p(k-1)}|\Delta
u(k-1)|^{p(k-1)}=+\infty.
$$
\end{lemma}
\begin{pf}
\noindent By \cite[Lemma 1, part a)]{MRT} there exist two positive constants $C_1$ and $C_2$ such that
$$\sum_{k=1}^{T+1}|\Delta u(k-1)|^{p(k-1)}\geq C_1\|u\|^{p^-}-C_2,$$
\noindent for every $u\in H$ with $\|u\|>1$. Hence, the conclusion simply holds taking into account that
\begin{eqnarray*}
 \Phi(u)&:=& \sum_{k=1}^{T+1}\frac{1}{p(k-1)}|\Delta
u(k-1)|^{p(k-1)}\\ &\geq &
\frac{1}{p^+}\left(\sum_{k=1}^{T+1}|\Delta u(k-1)|^{p(k-1)}\right)\\ &\geq&
 \frac{C_1\|u\|^{p^-}-C_2}{p^+}\rightarrow +\infty,
\end{eqnarray*}
as $\|u\|\rightarrow \infty$.
\end{pf}

 \indent A \textit{critical point} for $J_\lambda$, i.e. such a point $u\in H$ such that
 $$
 \langle J_\lambda^{'}(u),v\rangle=0,
 $$
 for every $v\in H$,
is said to be a \textit{weak solution} to \eqref{ani}.\par
\indent Summing by parts we see that any {weak solution} to \eqref{ani} is in fact a strong one. Hence, in order to
solve \eqref{ani} we check to find critical points for $J_\lambda$ and investigate their multiplicity.\par
 Finally, let us fix a constant $p>1$. Then
\begin{equation}\label{immersione}
\|u\|_\infty:=\max_{k\in \ZZ[1,T]}|u(k)|\leq \frac{(T+1)^{(p-1)/p}}{2}\left(\sum_{k=1}^{T+1}|\Delta
u(k-1)|^p\right)^{1/p},
\end{equation}
for every $u\in H$.\par
 See Lemma 2.2 proved by Jiang and Zhou in \cite{JiangZhou2}.
\section{Main Results}

Set
\[
\displaystyle{A_{0}:=\liminf_{t\rightarrow 0^+}}\frac{\displaystyle{\sum_{k=1}^{T}}\displaystyle{\max_{|\xi|\leq
t}F_k(\xi)}}{t^{p^+}}, \quad\quad\quad  \quad B^{0}:=\limsup_{t\rightarrow 0^+}\displaystyle{\frac{\displaystyle\sum_{k=1}^{T}F_k(t)}{t^{p^-}}},
\]
and
$$
\kappa:=\displaystyle\frac{2^{p^+-1}p^-}{p^+(T+1)^{p^+-1}}.
$$

\indent Our main result is the following one.
\begin{theorem}\label{Ma}
Assume that the following inequality
holds
\begin{itemize}
    \item [$(\textrm{h}_0)$] \qquad \qquad $A_0 <\kappa B^0.$
\end{itemize}

\noindent Then, for each $\lambda \in \displaystyle\left]
\frac{2}{p^-B^0},
\frac{2^{p^+}}{p^+(T+1)^{p^+-1}A_0}\right[$, problem
\eqref{ani} admits a sequence of nonzero solutions which converges to zero.
\end{theorem}

\begin{pf}
Fix $\lambda \in \displaystyle\left]
\frac{2}{p^-B^0},
\frac{2^{p^+}}{p^+(T+1)^{p^+-1}A_0}\right[$ and put
$\Phi$, $\Psi$,  $J_\lambda$ as
in the previous section.
Our aim is to apply Theorem \ref{two} to function $J_\lambda$.
Now, by standard arguments and bearing in mind Lemma \ref{LL}, condition $(\Lambda)$ clearly holds. Therefore,  our conclusion 
follows provided that $\delta<+\infty$  as well as that $0_H$ is not a local minimum point for the functional $J_\lambda$. To this end, let $\{c_m\}\subset ]0,+\infty[$ be a sequence such
that $\displaystyle \lim_{m\rightarrow \infty}c_m=0$ and
\[
\displaystyle{\lim_{m\rightarrow \infty}}\frac{\displaystyle{\sum_{k=1}^{n}}\displaystyle{\max_{|\xi|\leq
c_m}F_k(\xi)}}{c_m^2}=A_0.
\]
Set \[r_m:=\left(\frac{2^{p^+}}{p^+(T+1)^{p^+-1}}\right)c_m^{p^+},\]
for every $m\in\enne$.\par
\noindent For $m>0$ sufficiently large, it follows that
\begin{equation}\label{inc}
\Phi^{-1}(]-\infty,r_m[)\subset \{v\in X:|v(k)|\leq c_m,\,\,\forall\, k\in \ZZ[0,T+1]\}.
\end{equation}
Indeed, if $v\in X$ and $\Phi(v)<r_m$, one has
$$
\sum_{k=1}^{T+1}\frac{1}{p(k-1)}|\Delta v(k-1)|^{p(k-1)}<r_m.
$$
Then,
$$
|\Delta v(k-1)|<(p(k-1)r_m)^{1/p(k-1)},
$$
for every $k\in \ZZ[1,T+1]$. Consequently, since $r_m<1/p^+$ for every $m\geq \bar{m}$, one immediately one has
$$
(p^{+}r_m)^{1/p(k-1)}\leq (p^{+}r_m)^{1/p^{+}},
$$
hence
$$
|\Delta v(k-1)|<(p^{+}r_m)^{1/p^{+}}<1,
$$
for every $k\in \ZZ[1,T+1]$. So
$$
\sum_{k=1}^{T+1}\frac{1}{p^+}|\Delta v(k-1)|^{p^{+}}\leq \sum_{k=1}^{T+1}\frac{1}{p(k-1)}|\Delta v(k-1)|^{p(k-1)}<r_m,
$$
for every $m\geq \bar{m}$.
At this point, by \eqref{immersione}, the following inequality
$$
\max_{k\in \ZZ[1,T]}|v(k)|\leq \left(p^+\frac{(T+1)^{p^+-1}}{2^{p^+}}r_m\right)^{1/p^+}=c_m,
$$
\noindent is satisfied for every $m\geq \bar m$. Thus, the above computations ensure that \eqref{inc} holds. The algebraic inclusion \eqref{inc} implies that
$$
\varphi(r_m)\leq \frac{\displaystyle\sup_{\Phi(v)<
r_m}\displaystyle\sum_{k=1}^{T}F_k(v(k))}{r_m}
\leq
\frac{\displaystyle\sum_{k=1}^{T}\max_{|t|\leq c_m}F_k(t)}{r_m}=\frac{p^+(T+1)^{p^+-1}}{2^{p^+}}\frac{\displaystyle\sum_{k=1}^{T}\max_{|t|\leq c_m}F_k(t)}{c_m^{p^+}},
$$
for every $m$ sufficiently large.
Hence, it follows that
$$
\delta\leq \lim_{m\rightarrow \infty}\varphi(r_m)\leq\displaystyle \frac{p^+(T+1)^{p^+-1}}{2^{p^+}} A_0 < \frac{1}{\lambda}<+\infty.
$$
\noindent In the next step we show that $0_H$ is not a local minimum point for the functional $J_{{\lambda}}$.
First, assume that $B^0=+\infty$. Accordingly, fix $M$ such that $M>\displaystyle \frac{2}{p^-\lambda}$ and let $\{b_m\}$ be
a sequence of positive numbers,
with $\displaystyle\lim_{m\rightarrow \infty}b_m=0$, such that
\[
\sum_{k=1}^{T}F_k(b_m)>Mb^{p^-}_m \quad \quad  (\forall\; m\in \enne).
\]
Thus, take in $H$ a sequence $\{s_m\}$ such that, for each $m\in \enne$, $\displaystyle{s_m(k):=b_m}$ for every $k\in \ZZ[1,T]$.
Observe that $\|s_m\|\rightarrow 0$ as $m\rightarrow \infty$. So, it follows that
$$
\sum_{k=1}^{T+1}\frac{1}{p(k-1)}|\Delta
s_m(k-1)|^{p(k-1)}\leq \frac{2b_m^{p^-}}{p^-},
$$
taking into account that
$$
\sum_{k=1}^{T+1}|\Delta
s_m(k-1)|^{p^-}=2b_m^{p^-}.
$$
Then, one has
$$J_\lambda(s_m)\leq
\frac{2b_m^{p^-}}{p^-} - \lambda\sum_{k=1}^{T}F_k(b_m)<\left(\frac{2}{p^-}-{\lambda}M \right)b^{p^-}_m,$$
that is, $J_\lambda(s_m)<0$ for every sufficiently large $m$.

\noindent Next, assume that $B^0<+\infty$. Since $\lambda >\displaystyle\frac{2}{p^-B^0}$, we can fix $\varepsilon>0$ such that
${\varepsilon<B^0-\displaystyle \frac{2}{p^-\lambda}}$. Therefore,
also taking $\{b_m\}$ a sequence of positive numbers such that $\displaystyle\lim_{m\rightarrow \infty}b_m=0$ and
$$
(B^0-\varepsilon)b^{p^-}_m
<\sum_{k=1}^{T}F_k(b_m)<(B^0+\varepsilon) b^{p^-}_m,  \quad \quad  (\forall\; m\in \enne)
$$
arguing as before and by choosing $\{s_m\}$ in $H$ as above, one has
$$J_\lambda(s_m)<\left(\frac{2}{p^-}-{\lambda}(B^0-\varepsilon) \right) b^{p^-}_m.$$
So, also in this case, $J_\lambda(s_m)<0$ for every sufficiently large $m$. Finally, since $J_{{\lambda}}(0_H)=0$, the above fact means that $0_H$ is not a local minimum of $J_{{\lambda}}$. Therefore, the unique global minimum of $\Phi$ is not a local minimum of the functional $J_{{\lambda}}$. Hence, by Theorem \ref{two} we obtain a sequence $\{u_m\}\subset H$ of critical points of $J_{{\lambda}}$ such that
$$\lim_{m\rightarrow\infty}\sum_{k=1}^{T+1}\frac{1}{p(k-1)}|\Delta
u_m(k-1)|^{p(k-1)}=\lim_{m\rightarrow\infty} \|u_m\|=0.$$ Thus, it follows that $\|u_m\|_\infty\rightarrow 0$ as $m\rightarrow \infty$. The proof is complete.
\end{pf}

\begin{remark}\label{intervallicommento}
\rm{As pointed out earlier, Theorem \ref{Ma} has been obtained exploiting Theorem \ref{two}. Via our approach we are able to determine an open subinterval of $]0,1/\delta[$, where $\delta:=\displaystyle\liminf_{r\rightarrow 0^+}\varphi(r)$, such that problem \eqref{ani} admits infinitely many solutions. Indeed, the main condition $(\textrm{h}_0)$ implies that the real interval of parameter
$$\displaystyle\left]
\frac{2}{p^-B^0},
\frac{2^{p^+}}{p^+(T+1)^{p^+-1}A_0}\right[,$$
is well-defined and nonempty. Further, since $B_0>0$ and
$$
\delta\leq \lim_{m\rightarrow \infty}\varphi(r_m)\leq\displaystyle \displaystyle \frac{p^+(T+1)^{p^+-1}}{2^{p^+}} A_0,
$$
the following inclusion
$$
\displaystyle\left]
\frac{2}{p^-B^0},
\frac{2^{p^+}}{p^+(T+1)^{p^+-1}A_0}\right[\subseteq \left]0,\frac{1}{\delta}\right[,
$$
is verified.
}
\end{remark}

\begin{remark}\label{cp}
\rm{We note that, if $f_k$ is a nonnegative continuous function, for every $k\in \ZZ[1,T]$, condition $(\textrm{h}_0)$ assumes the form
\[
{\liminf_{t\rightarrow 0^+}}\frac{\displaystyle \sum_{k=1}^{n}F_k(t)}{t^{p^+}}<\kappa\limsup_{t\rightarrow 0^+}\frac{\displaystyle \sum_{k=1}^{n}F_k(t)}{t^{p^-}}.
\]
Consequently, Theorem \ref{intro} immediately follows from Theorem \ref{Ma}.
}
\end{remark}

We note that there seems to be increasing interest in existence of solutions to boundary value
problems for finite difference equations with $p$-Laplacian operator, because of their applications
in many fields. In this setting, set $p>1$ and consider the real map $\phi_p:\erre\rightarrow \erre$ given by $\phi_p(s):=|s|^{p-2}s$, for every $s\in\erre$.\par Further, denote
\[
\displaystyle{\widehat{A}_{0}:=\liminf_{t\rightarrow 0^+}}\frac{\displaystyle{\sum_{k=1}^{T}}\displaystyle{\max_{|\xi|\leq
t}F_k(\xi)}}{t^{p}} \quad  \quad\textrm{ and} \quad  \quad \widehat{B}^{0}:=\limsup_{t\rightarrow 0^+}\displaystyle{\frac{\displaystyle\sum_{k=1}^{T}F_k(t)}{t^{p}}}.
\]
With the previous notations, taking the map $p:\ZZ[0,T+1]\rightarrow \erre$ such that $p(k)=p,$ for every $k\in \ZZ[0,T+1],$ we have the following immediate consequence of Theorem \ref{Ma}.
\begin{corollary}\label{Ma3}
Assume that
\begin{itemize}
    \item [$(\widehat{\textrm{h}}_0)$] \qquad \qquad $\widehat{A}_0 <\displaystyle \frac{2^{p-1}}{(T+1)^{p-1}} \widehat{B}^0.$
\end{itemize}

\noindent Then, for each $\lambda \in \displaystyle\left]
\frac{2}{p\widehat{B}^0},
\frac{2^{p}}{p(T+1)^{p-1}\widehat{A}_0}\right[$, the following problem
\begin{equation}\tag{$D_{\lambda}^{f}$} \label{dani3}
\left\{
\begin{array}{l}
-\Delta(\phi_p(\Delta u({k-1})))
= \lambda f_k(u(k)) ,\quad k \in
\ZZ[1,T],
\\ {u(0)=u({T+1})=0, } \quad \quad \qquad \qquad \qquad \\
\end{array}
\right.
\end{equation}
\noindent admits a sequence of nonzero solutions which converges to zero.
\end{corollary}

A more technical version of Theorem \ref{Ma} can be written as follows.
\begin{theorem}\label{main2}
Assume that there exist nonnegative real sequences $\{a_m\}$ and $\{b_m\}$, with $\displaystyle\lim_{m\rightarrow \infty}b_m=0$, such that
\begin{itemize}
    \item [$(\textrm{k}_1)$] \quad ${a_m^{p^-}}<\displaystyle\left(\frac{2^{p^+-1}p^-}{p^+(T+1)^{p^+-1}}\right)b_m^{p^+},$
     for each $ m\in \enne$$;$
     \item [$(\textrm{k}_2)$]\quad $\displaystyle G_0<
   {\frac{B^0}{p^+(T+1)^{p^+-1}}}$, where
   $$
   G_0:=\lim_{m\rightarrow \infty}\frac{\displaystyle{\sum_{k=1}^{T}\max_{|t|\leq b_m}F_k(t)}-\sum_{k=1}^{T}F_k(a_m)}{\displaystyle 2^{p^+-1}p^-b_m^{p^+}-p^+(T+1)^{p^+-1}a_m^{p^-}}.
   $$
\end{itemize}
Then, for each $\lambda \in \displaystyle\left] \frac{2}{p^-B^0},
\frac{2}{p^-p^+(T+1)^{p^+-1}G_0}\right[$, problem
\eqref{ani} admits a sequence of nonzero solutions which converges to zero.
\end{theorem}
\begin{pf}Let us observe that, keeping the above notations, one has
\begin{equation}\label{jfk}
\varphi(r_m)\leq \inf_{w \in \Phi^{-1}(\left]-\infty, r_m\right[)}\frac{\displaystyle{\sum_{k=1}^{T}\max_{|t|\leq b_m}F_k(t)}-\sum_{k=1}^{T}F_k(w(k))}{\displaystyle\left(\frac{2^{p^+}}{p^+(T+1)^{p^+-1}}\right)b_m^{p^+}-\Phi(u)}.
\end{equation}
Now,
for each $m\in \enne$, let
$w_m\in H$
be defined by $w_m(k):=a_m$ for every $k\in \ZZ[1,T]$.
Clearly, since $\|w_m\|\rightarrow 0$ as $m\rightarrow \infty$, it follows that $|\Delta w_m(k-1)|<1$ for every $k\in \ZZ[1,T+1]$ and sufficiently large $m$. Then, there exists $\bar m\in\enne$ such that
$$
\Phi(w_m):=\sum_{k=1}^{T+1}\frac{1}{p(k-1)}|\Delta
w_m(k-1)|^{p(k-1)}\leq \frac{1}{p^-}\left(\sum_{k=1}^{T+1}|\Delta
w_m(k-1)|^{p^-}\right),
$$
for every $m\geq \bar m$.\par
\noindent Now, by condition $(\textrm{k}_1)$ and taking into account that $\displaystyle\sum_{k=1}^{T+1}|\Delta w_m(k-1)|^{p^{-}}=2a_m^{p^-}$, the above inequality implies that
$$
0<\displaystyle\left(\frac{2^{p^+}}{p^+(T+1)^{p^+-1}}\right)b_m^{p^+}-\frac{2a_m^{p^-}}{p^-}\leq \displaystyle\left(\frac{2^{p^+}}{p^+(T+1)^{p^+-1}}\right)b_m^{p^+}-\Phi(w_m),
$$
\noindent for every $m\geq \bar m$.\par
\noindent Hence, clearly $w_m \in \Phi^{-1}(\left]-\infty, r_m\right[)$ and inequality \eqref{jfk} yields
\[
\varphi(r_m)\leq \left(\frac{p^-p^+(T+1)^{p^+-1}}{2}\right)\frac{\displaystyle{\sum_{k=1}^{T}\max_{|t|\leq b_m}F_k(t)}-\sum_{k=1}^{T}F_k(a_m)}{\displaystyle p^-2^{p^+-1}b_m^{p^+}-p^+(T+1)^{p^+-1}a_m^{p^-}},
\]
\noindent for every large enough $m$.\par
\noindent Further, by hypothesis $(\textrm{k}_2)$, we obtain
$$
\delta\leq \lim_{m\rightarrow \infty}\varphi(r_m)\leq\displaystyle \frac{p^-p^+(T+1)^{p^+-1}}{2}G_0 < \frac{1}{\lambda}<+\infty.
$$
From now on, arguing exactly as in the proof of Theorem \ref{Ma} we obtain the assertion.
\end{pf}

 A more general form of our results can be obtained by using a lemma due to Cabada, Iannizzotto and Tersian; see \cite[Lemma 4]{cit}. Indeed, in the cited work, the authors improved \cite[Lemma 2.2]{JiangZhou2}, finding that
\begin{equation}\label{immersione2}
\|u\|_\infty\leq \frac{1}{c_1}\left(\sum_{k=1}^{T+1}|\Delta
u(k-1)|^p\right)^{1/p},
\end{equation}
for every $u\in H$, where
\begin{equation*}
c_1:=\left\{
\begin{array}{ll}
\left[\displaystyle\left(\frac{2}{T}\right)^{p-1}+\left(\frac{2}{T+2}\right)^{p-1}\right]^{1/p} & {\rm if}\,\,\, T\,\, {\rm is\,\,even} \\
 &\\
\displaystyle \frac{2}{(T+1)^{(p-1)/p}} & {\rm if}\,\,\, T\,\, {\rm is\,\,odd.} \\
\end{array}
\right.
\end{equation*}
Note that, since the continuous function $\theta:]0,T+1[\rightarrow \erre$ defined by
$$
\theta(s):=\frac{1}{(T-s+1)^{p-1}}+\frac{1}{s^{p-1}},
$$
attain its minimum $\displaystyle\frac{2^p}{(T+1)^{p-1}}$ at $\displaystyle s=\frac{T+1}{2}$, one has
$$
\displaystyle\frac{2^p}{(T+1)^{p-1}}< \theta(T/2).
$$
Hence
\begin{equation}\label{immersione3}
\frac{2}{(T+1)^{(p-1)/p}}< \left[\displaystyle\left(\frac{2}{T}\right)^{p-1}+\left(\frac{2}{T+2}\right)^{p-1}\right]^{1/p}=\theta(T/2)^{1/p}.
\end{equation}
Then
$$
\left[\left(\frac{2}{T}\right)^{p-1}+\left(\frac{2}{T+2}\right)^{p-1}\right]^{-1/p}<\frac{(T+1)^{(p-1)/p}}{2}.
$$

For instance, an immediate consequence of the previous remarks is the following result.

\begin{corollary}\label{Ma5} Let $T\geq 2$ be an even even number.
Assume that
\begin{itemize}
    \item [$(\widetilde{{\textrm{h}}}_0)$] \qquad \qquad $\widehat{A}_0 <\displaystyle 2^{p-2}\left[\displaystyle\frac{1}{T^{p-1}}+\frac{1}{(T+2)^{p-1}}\right] \widehat{B}^0.$
\end{itemize}

\noindent Then, for each $\lambda \in \displaystyle\left]
\frac{2}{p\widehat{B}^0},
\left[\displaystyle\left(\frac{2}{T}\right)^{p-1}+\left(\frac{2}{T+2}\right)^{p-1}\right]\frac{1}{p\widehat{A}_0}\right[$, the following problem
\begin{equation}\tag{$D_{\lambda}^{f}$} \label{dani2}
\left\{
\begin{array}{l}
-\Delta(\phi_p(\Delta u({k-1})))
= \lambda f_k(u(k)) ,\quad k \in
\ZZ[1,T]
\\ {u(0)=u({T+1})=0, } \quad \quad \qquad \qquad \qquad \\
\end{array}
\right.
\end{equation}
\noindent admits a sequence of nonzero solutions which converges to zero.
\end{corollary}

\begin{remark}\label{Confronto}
\rm{It follows from \eqref{immersione3} that condition $(\widehat{\textrm{h}}_0)$ imply $(\widetilde{{\textrm{h}}}_0)$ and one has
$$
\displaystyle\left]
\frac{2}{p\widehat{B}^0},
\frac{2^{p}}{p(T+1)^{p-1}\widehat{A}_0}\right[\subseteq \left]
\frac{2}{p\widehat{B}^0},
\left[\displaystyle\left(\frac{2}{T}\right)^{p-1}+\left(\frac{2}{T+2}\right)^{p-1}\right]\frac{1}{p\widehat{A}_0}\right[.
$$
Moreover, if $\widehat{A}_0>0$, again inequality \eqref{immersione3} ensures that the above inclusion is proper.
Hence, from the previous facts, we can conclude that Corollary \ref{Ma5} is a refinement of Corollary \ref{Ma3}.
}
\end{remark}
An application of Theorem \ref{intro} is the following.
\begin{example}\label{esempio}
{\rm Let $\gamma>2$ be a real positive constant and $p:\ZZ[0,T+1]\rightarrow \erre^+$ a map such that $p^+:=\displaystyle \max_{\ZZ[0,T+1]}p(k)=\gamma$ and $p^-:=\displaystyle \min_{\ZZ[0,T+1]}p(k)=\gamma-1$, with $T\geq 2$.
Further, let
$\{s_m\}$, $\{t_m\}$ and $\{\delta_m\}$ be sequences defined by
$$
\displaystyle s_m:=2^{-\frac{m!}{2}},\,\,\,\,\,\,\,\,\
t_m:=2^{-2m!},\,\,\,\,\,\,\,\, \delta_m:=2^{-(m!)^2},
$$
\noindent and consider $\nu\in \enne$ such that
$$
s_{m+1}<t_m<s_m-\delta_m,\,\,\,\,\,\,\forall\; m\geq \nu.
$$
Moreover, let $g:\RR\rightarrow\RR$ be the
nonnegative continuous function given by
\begin{equation*}
g(s):=\left\{
\begin{array}{ll}
2^{-(\gamma-1)\nu!}& {\rm if} \quad s\in ]s_{\nu}-\delta_\nu,+\infty[ \\
y_m(s) &{\rm if} \quad s\in \bigcup_{m\geq \nu}]s_{m+1}-\delta_{m+1},s_{m+1}[\\
2^{-(\gamma-1)m!}& {\rm if} \quad s\in \bigcup_{m\geq \nu}[s_{m+1},s_{m}-\delta_m] \\
0 & {\rm if}\quad  s\leq 0,\\
\end{array}
\right.
\end{equation*}
\noindent where
$$
y_m(s):=\displaystyle\left(2^{-(\gamma-1)m!}-2^{-(\gamma-1)(m+1)!}\right)\left(\frac{s-s_{m+1}+\delta_{m+1}}{\delta_{m+1}}\right)+2^{-(\gamma-1)(m+1)!}.
$$
 \noindent Set  $\displaystyle
G(t):=\int_0^t g(s)\;ds$ for every $t\in\erre$. Then, one has that
$$
\frac{G(s_m)}{s_m^{\gamma}}\leq
\frac{g(s_{m+1})s_m+g(s_m)\delta_m}{s_m^{\gamma}},
$$
\noindent and
$$
\frac{G(t_m)}{t_m^{\gamma-1}}\geq
\frac{g(s_{m+1})(t_m-s_{m+1})}{t_m^{\gamma-1}},
$$
\noindent for every large enough $m$.\par
\noindent Since
$$
\lim_{m\rightarrow\infty}\frac{g(s_{m+1})s_m+g(s_m)\delta_m}{s_m^{\gamma}}=0,\quad\textrm{and}\quad
\lim_{m\rightarrow\infty}\frac{g(s_{m+1})(t_m-s_{m+1})}{t_m^{\gamma-1}}=+\infty,
$$
\noindent it follows that
$$
\lim_{m\rightarrow\infty}\frac{G(s_m)}{s_m^\gamma}=0, \quad\textrm{and}\quad
\lim_{m\rightarrow\infty}\frac{G(t_m)}{t_m^{\gamma-1}}=+\infty.
$$
Thus
\[
\displaystyle{\liminf_{t\rightarrow 0^+}}\frac{\displaystyle\int_{0}^{t}g(s)ds}{t^{\gamma}}=0,\quad\textrm{and}\quad\limsup_{t\rightarrow 0^+}\frac{\displaystyle\int_{0}^{t}g(s)ds}{t^{\gamma-1}}=+\infty.
\]
Then, owing to Theorem \ref{intro}, for each $\lambda>0$, the following anisotropic discrete Dirichlet problem
\begin{equation}\tag{$A_{\lambda}^{g}$} \label{ani3}
\left\{\begin{array}{lll}
-\Delta(|\Delta u(k-1)|^{p(k-1)-2}\Delta u(k-1))=\lambda g(u(k)),\,\,
k\in\ZZ[1,T]\\
u(0)=u(T+1)=0,
\end{array}\right.
\end{equation}
admits a sequence of nonzero solutions which converges to zero.}
\end{example}

 \medskip
 \indent {\bf Acknowledgements.} The research was supported in part by the SRA grants P1-0292-0101 and J1-4144-0101. The paper was realized with the auspices of the GNAMPA Project 2012 entitled: {\it Esistenza e molteplicit\`{a} di soluzioni per problemi differenziali non lineari}.

\end{document}